\documentclass[12 pt]{amsart}
\usepackage{amssymb}
\newtheorem{theorem}{Theorem}[]
\newtheorem{lemma}[theorem]{Lemma}
\newtheorem{proposition}[theorem]{Proposition}
\newcommand{\mH}{\mathcal{H}}
\newcommand{\M}{\mathcal{M}}
\newcommand{\e}{\ell^1}
\newcommand{\gd}{\delta}
\def\t{\times}
\def\ot{\otimes}
\def\ho{\hat{\otimes}}
\def\d{\cdot}
\def\to{\longrightarrow}
\begin{document}
\baselineskip15pt
\title[Morita Equivalence]{Morita Equivalence of Brandt Semigroup Algebras}
\author[M.M. Sadr]{Maysam Maysami Sadr}
\address{Department of Mathematics, Institute for Advanced Studies in Basic Sciences, Zanjan, Iran}
\email{sadr@iasbs.ac.ir}
\begin{abstract}
We prove that for every group $G$ and any two sets $I,J$,
the Brandt semigroup algebras $\e(B(I,G))$ and $\e(B(J,G))$ are Morita equivalent
with respect to the Morita theory of self-induced Banach algebras introduced by Gr{\o}nb{\ae}k.
As applications, we show that if $G$ is an amenable group, then for a wide class of Banach
$\e(B(I,G))$-bimodules $E$, and every $n>0$, the bounded Hochschild cohomology groups
$\mH^n(\e(B(I,G)),E^*)$ are trivial, and also, the notion of approximate amenability
is not Morita invariant.
\end{abstract}
\maketitle
\section{Introduction}
Morita theory is a very useful tool in the study of rings and algebras. In the area of topological algebras,
there are different notions of Morita theory in the literature, but all of them has no many variations
with respects to the original theory, defined by Kiiti Morita, in pure algebraic case.
Niels Gr{\o}nb{\ae}k defined a Morita theory for Banach algebras with bounded approximate
identities in \cite{G1}. Then he extended his theory to the larger class of self-induced Banach algebras
in \cite{G2}. In the first theory, the only class of algebras that are Morita equivalent to the algebra of
complex numbers $\mathbb{C}$, is the class of finite dimensional matrix algebras, but in the
second one we find many infinite dimensional  Banach algebras that are Morita equivalent to $\mathbb{C}$.
In this paper we construct a new class of infinite dimensional self-induced Banach algebras Morita
equivalent to $\mathbb{C}$. Then by this construction, we show that for a discrete group $G$, and every two
nonempty sets $I,J$, the Banach convolution algebras $\e(B(I,G))$ and $\e(B(J,G))$ are Morita equivalent,
where $B(I,G)$ denotes the Brandt semigroup over $G$ with index set $I$, \cite{Ho}.
Brandt semigroups are one of the most
important inverse semigroups. Some authors have payed attention to the bounded Hochschild cohomology and amenability properties of inverse semigroup algebras, \cite{DN}, \cite{DP}, \cite{G5}, \cite{C}, \cite{GLZ}, \cite{LS}, \cite{MP} and \cite{M1}. As a corollary of Morita equivalence of Brandt semigroup algebras and a strong result of \cite{G2},
we show that if $G$ is an amenable semigroup and $I$ is an arbitrary nonempty set, then the topological Hochschild
homology and cohomology groups $\mH_n(\e(B(I,G)),E)$ and $\mH^n(\e(B(I,G)),E^*)$, for any $n>0$ and every
$\e(B(I,G))$-induced Banach bimodule $E$, are trivial. Also by a specific example, we show that
the notion of approximate amenability of Banach algebras, introduced by Ghahramani and Loy \cite{GL},
is not Morita invariant. This result is against Gr{\o}nb{\ae}k's corollary on amenability \cite[Corollary 6.5]{G1}
that says for Banach algebras with bounded approximate identities the notion of amenability is
Morita invariant.
\section{Preliminaries}
Throughout for an element $x$ of a set $X$, $\gd_x$ is its
point mass measure in $\e(X)$. Let $E$ and $F$ be Banach spaces.
The Banach space which is the completed projective tensor product of
$E$ and $F$ is denoted by $E\ho F$; for $z\in E\ho F$, there are sequences
$(x_n)\in E$ and $(y_n)\in F$ such that $\sum_{n=1}^\infty\|x_n\|\|y_n\|<\infty$ and
$z=\sum_{n=1}^\infty x_n\ot y_n$. Similar pure algebraic case, the key property of
$E\ho F$ is that, for each continuous bilinear map $B:E\t F\to D$, where $D$ is a
Banach space, there is a unique continuous linear map $T:E\ho F\to D$ with $\|T\|=\|B\|$
and $T(x\ot y)=B(x,y)$ ($x\in E,y\in F$), see \cite{H} for more details.

It is well known that for non empty sets $X$ and $Y$, the map $\gd_x\ot \gd_y\longmapsto\gd_{(x,y)}$
($x\in X, y\in Y$), defines an isometric isomorphism between Banach spaces $\e(X)\ho\e(Y)$ and
$\e(X\t Y)$; we use frequently this identification.

Let $A, B$ and $C$ be Banach algebras. A left Banach $A$-module $E$, is an ordinary
left $A$-module which is a Banach space and there is a constant $M>0$ such that
$\|a\d x\|\leq M\|a\| \|x\|$ ($a\in A,x\in E$). Similarly , right Banach modules
and Banach bimodules are defined. The category of left Banach $A$-modules
and bounded module homomorphisms is denoted by $A$-$\mathbf{mod}$.
Similarly, one can define the category of right Banach $A$-modules
$\mathbf{mod}$-$A$, and the category of left-$A$ right-$B$ Banach bimodules
$A$-$\mathbf{mod}$-$B$. The notations $_AE$, $E_A$ and $_AE_B$ are shorthand
indications that $E$ is in $A$-$\mathbf{mod}$, $\mathbf{mod}$-$A$ and $A$-$\mathbf{mod}$-$B$, respectively.

For $E$ in $\mathbf{mod}$-$A$ and $F$ in $A$-$\mathbf{mod}$, let $E\ho_A F$ be the universal object for
$A$-balanced bounded bilinear maps from $E\t F$. This can be realized as the Banach space
$E\ho F/N$ where $N$ is the closed linear span of $\{x\d a\ot y-x\ot a\d y: x\in E, y\in F, a\in A\}$
(\cite{R}, \cite{H}). For $_BE_A$ and $_AF_C$ the tensor product $E\ho_A F$ is in $B$-$\mathbf{mod}$-$C$.
For $_AE$, define a left module homomorphism $\mu_E:A\ho_A E\to E$ by $\mu_E(a\ot x)=a\d x$
($a\in A,x\in E$).

The Banach algebra $A$ is called self-induced if the multiplication map $a\ot b\to ab$,
from $A\ho_A A$ to $A$, is an isomorphism (between $A$-bimodules). More generally, a left Banach
$A$-module $E$, is $A$-induced if $\mu_E$ is an isomorphism (\cite{G3}). Similarly, right and
two sided induced modules are defined. Denote by ind-$A$-$\mathbf{mod}$, the full subcategory
of $A$-$\mathbf{mod}$, of left Banach $A$-induced modules.

Two self-induced Banach algebras $A$ and $B$ are called Morita equivalent, and denoted by
$A\approx B$, if there are covariant functors
$$\Phi:\text{ind-}A\text{-}\mathbf{mod}\to\text{ind-}B\text{-}\mathbf{mod}\quad
\Psi:\text{ind-}B\text{-}\mathbf{mod}\to\text{ind-}A\text{-}\mathbf{mod}$$
with appropriate continuity properties such that $\Psi\Phi\simeq id_{\text{ind-}A\text{-}\mathbf{mod}}$
and $\Phi\Psi\simeq id_{\text{ind-}B\text{-}\mathbf{mod}}$, where $\simeq$ denotes natural isomorphism
between functors, for complete definition see the original paper \cite{G2}.
We only need the following characterization of Morita equivalent Banach algebras \cite{G2}.
\begin{theorem}\label{t1}
Let $A$ and $B$ be self-induced Banach algebras. Then $A$ and $B$ are Morita equivalent
if and only if there are two sided induced modules $P\in B$-$\mathbf{mod}$-$A$ and
$Q\in A$-$\mathbf{mod}$-$B$ such that $P\ho_A Q\cong B$ and $Q\ho_B P\cong A$, where $\cong$
denotes topological isomorphism of bimodules.
\end{theorem}
\section{A  Banach algebra Morita equivalent to $\mathbb{C}$}
In this section, for any set $I$, we define a {\it matrix-like} Banach algebra $\M_I$
and prove that $\M_I$ is Morita equivalent to the algebra of complex numbers $\mathbb{C}$.

Let $I$ be a nonempty set. Let the underlying Banach space of $\M_I$ be  $\e(I\times I)$
and let its multiplication be the convolution product
$$(ab)(i,j)=\sum_{k\in I}a(i,k)b(k,j)\quad(a,b\in\M_I, i,j\in I).$$
Note that if $I$ is a finite set, then $\M_I$ is isomorphic to an ordinary matrix algebra.
Also, for any $i,p,q,j\in I$, we have the following identity in $\M_I$:
$$\gd_{(i,p)}\gd_{(q,j)}=\begin{cases}
\gd_{(i,j)} &\hbox{\rm\ if }p=q\\
0&\hbox{\rm\ if }p\neq q.\\
\end{cases}$$
Define a two sided Banach module action of $\M_I$ on $\e(I)$ by
$$(a\d b)(i)=\sum_{k\in I}a(i,k)b(k),\quad(b\d a)(i)=\sum_{k\in I}b(k)a(k,i),$$
for $a\in\M_I,b\in\e(I),i\in I$.
\begin{lemma}\label{l1}
The map $\nu:\e(I)\ho_{\M_I}\e(I)\to\mathbb{C}$, defined by $\nu(a\ot b)=\sum_{i\in I}a(i)b(i)$,
is an isomorphism of Banach spaces.
\end{lemma}
\begin{proof}
By definition, it is enough to prove that the map $$\bar{\nu}:\e(I)\ho\e(I)\to\mathbb{C},$$
defined by $\bar{\nu}(a\ot b)=\sum_{i\in I}a(i)b(i)$, is nonzero and $N$, the closed linear span of
$\{\gd_i\d\gd_{(j,k)}\ot\gd_{i'}-\gd_i\ot\gd_{(j,k)}\d\gd_{i'}:i,i',j,k\in I\}$, is equal to $ker\bar{\nu}$.\\
If $k_0\in I$ is arbitrary, then $\bar{\nu}(\gd_{k_0}\ot\gd_{k_0})=1$. This shows that $\bar{\nu}$ is not zero.\\
A simple computation shows that for every $i,i',j,k\in I$, $\bar{\nu}(\gd_i\d\gd_{(j,k)}\ot\gd_{i'})=
\bar{\nu}(\gd_i\ot\gd_{(j,k)}\d\gd_{i'})$. This implies that $N\subseteq ker\bar{\nu}$.\\
For the converse, we have
\begin{equation}\label{e1}
\gd_{(i,j)}\in N,\quad\quad\text{if }\quad i\neq j,
\end{equation}
since $\gd_i\d\gd_{(j,k)}\ot\gd_{k}-\gd_i\ot\gd_{(j,k)}\d\gd_k=-\gd_i\ot\gd_j$. Also for every $i,j\in I$,
we have
\begin{equation}\label{e2}
\gd_{(j,j)}-\gd_{(i,i)}\in N,
\end{equation}
since $\gd_i\d\gd_{(i,j)}\ot\gd_{j}-\gd_i\ot\gd_{(i,j)}\d\gd_j=\gd_j\ot\gd_j-\gd_i\ot\gd_i$.
Now suppose that $c=\sum_{i,j\in I}c(i,j)\gd_{(i,j)}$ is in $ker\bar{\nu}$. Thus we have
\begin{equation}\label{e3}
\bar{\nu}(c)=\sum_{i\in I}c(i,i)=0.
\end{equation}
Consider the following decomposition of $c$ :
$$c=\sum_{i,j\in I, i\neq j}c(i,j)\gd_{(i,j)}+\sum_{i\in I}c(i,i)\gd_{(i,i)}=a+b.$$
Then by (\ref{e1}), $a$ is in $N$. Let $k_0\in I$ be arbitrary and fixed, then by (\ref{e3}),
we have $b=\sum_{i\in I}c(i,i)\gd_{(i,i)}-\sum_{i\in I}c(i,i)\gd_{(k_0,k_0)}$. Thus by (\ref{e2}),
$b$ is also in $N$. Therefore $c$ is in $N$ and $ker\bar{\nu}\subseteq N$.
\end{proof}
\begin{proposition}\label{p2}
\hspace{1mm}
\begin{enumerate}
\item [(i)] $\e(I)$ is a two sided $\M_I$-induced module.
\item [(ii)] $\M_I$ is a self-induced Banach algebra.
\end{enumerate}
\end{proposition}
\begin{proof}
The canonical map $\gd_i\ot \gd_j\longmapsto\gd_{(i,j)}$ ($i,j\in I$), from $\e(I)\ho\e(I)$ to $\M_I$,
is an isomorphism of Banach $\M_I$-bimodules. Thus we have,
\begin{equation*}
\begin{split}
\M_I\ho_{\M_I}\e(I)&\cong(\e(I)\ho\e(I))\ho_{\M_I}\e(I)\\
&\cong\e(I)\ho(\e(I)\ho_{\M_I}\e(I))\\
&\cong\e(I)\ho\mathbb{C}\hspace{40mm}\text{(by Lemma \ref{l1})}\\
&\cong\e(I).
\end{split}
\end{equation*}
This proves $\e(I)$ is left $\M_I$-induced. Similarly , it is proved that $\e(I)$ is
right $\M_I$-induced.
For (ii), we have,
\begin{equation*}
\begin{split}
\M_I\ho_{\M_I}\M_I&\cong(\e(I)\ho\e(I))\ho_{\M_I}\M_I\\
&\cong\e(I)\ho(\e(I)\ho_{\M_I}\M_I)\\
&\cong\e(I)\ho\e(I)\hspace{40mm}\text{(by (i))}\\
&\cong\M_I.
\end{split}
\end{equation*}
Thus $\M_I$ is self-induced.
\end{proof}
\begin{theorem}\label{t2}
$\M_I$ is Morita equivalent to $\mathbb{C}$.
\end{theorem}
\begin{proof}
By Lemma \ref{l1} and Proposition \ref{p2}, the Banach algebras $A=\M_I$, $B=\mathbb{C}$
and Banach bimodules $P=_\mathbb{C}\e(I)_{\M_I}$, $Q=_{\M_I}\e(I)_\mathbb{C}$ satisfy conditions
of Theorem \ref{t1}. Thus $\M_I$ is Morita equivalent to $\mathbb{C}$.
\end{proof}
{\bf Remark.}\begin{enumerate} \item[I.] It is proved in \cite{G2}, that for any Banach space $E$, the
{\it tensor algebra} $E\ho E^*$ is Morita equivalent to $\mathbb{C}$. Also, it is well known that if $E$ has
{\it bounded approximate property}, then the algebra $\mathcal{N}(E)$ of nuclear operators
on $E$, and $E\ho E^*$ are isomorph. Thus by Theorem \ref{t2}, $\mathcal{N}(\e(I))\cong\e(I)\ho\ell^\infty(I)$
and $\M_I$ are Morita equivalent, but clearly these are not isomorph if $I$ is an infinite set.
\item[II.] For Morita theory of some other Matrix-like algebras, see \cite{G1}, \cite{G2} and \cite{G4}.
\end{enumerate}
\section{The main result}
Let $I$ be a nonempty set and let $G$ be a discrete group. Consider the set $T=I\times G\times I$,
add an extra element $\o$ to $T$, and define a semigroup multiplication on $S=T\cup\{\o\}$,
as follows. For $i,i',j,j'\in I$ and $g,g'\in G$, let
$$(i,g,j)(i',g',j')=\begin{cases}
(i,gg',j') &\hbox{\rm\ if }j=i'\\
\o&\hbox{\rm\ if }j\neq i', \\
\end{cases}$$
also let $\o(i,g,j)=(i,g,j)\o=\o$ and $\o\o=\o$. Then $S$ becomes a semigroup
that is called Brandt semigroup over $G$ with index set $I$, and usually denoted by
$B(I,G)$. For more details see \cite{Ho}.\\
The Banach space $\e(T)$, with the convolution product
$$(ab)(i,g,j)=\sum_{k\in I,h\in G}a(i,gh^{-1},k)b(k,h,j),$$
for $a,b\in\e(T), i,j\in I, g\in G$ becomes a Banach algebra.
\begin{lemma}\label{l2}
The Banach algebras  $\e(S)$ and $\e(T)\oplus\mathbb{C}$ are homeomorphic isomorph,
where the multiplication of $\e(T)\oplus\mathbb{C}$ is coordinatewise.
\end{lemma}
\begin{proof}
Consider the following short exact sequence of Banach algebras and continuous algebra homomorphisms.
$$0\to\e(T)\to\e(S)\to\mathbb{C}\to0,$$
where the second arrow $u:\e(T)\to\e(S)$ is defined by $u(b)(t)=b(t)$ and $u(b)(\o)=-\sum_{s\in T}b(s)$,
for $b\in\e(T)$ and $t\in T\subset S$, and the third arrow $v:\e(S)\to\mathbb{C}$ is the
integral functional, $v(a)=\sum_{s\in S}a(s)$ ($a\in\e(S)$). Now, let $w:\e(S)\to\e(T)$ be the restriction
map, $w(a):=a|_T$. Then $w$ is a continuous  algebra homomorphism and $wu=Id_{\e(T)}$.
Thus the exact sequence splits and we have $\e(S)\cong\e(T)\oplus\mathbb{C}$.
\end{proof}
The following Theorem is our main result.
\begin{theorem}\label{t4}
Let $I$ and $J$ be nonempty sets and let $G$ be a discrete group. Then $\e(B(I,G))$ and $\e(B(J,G))$
are Morita equivalent self-induced Banach algebras.
\end{theorem}
\begin{proof}
Let $T$  be as above. It is easily checked that the map $\gd_{(i,j)}\ot\gd_g\longmapsto\gd_{(i,g,j)}$
($i,j\in I, g\in G$), is an isometric isomorphism from the Banach algebra $\M_I\ho\e(G)$ onto $\e(T)$. Thus
$\e(T)$ is self-induced, since $\M_I$ and $\e(G)$ are self-induced. Also, since $\M_I\approx\mathbb{C}$, we have
$\e(T)\approx\mathbb{C}\ho\e(G)\approx\e(G)$. By Lemma \ref{l2}, we have $\e(B(I,G))\cong\e(T)\oplus\mathbb{C}$
thus $\e(B(I,G))$ is self-induced and Morita equivalent to $\e(G)\oplus\mathbb{C}$. Similarly
$\e(B(J,G))\approx\e(G)\oplus\mathbb{C}$, therefore we have $\e(B(I,G))\approx\e(B(J,G))$.
\end{proof}
\section{Some applications}
For the topological Hochschild (co)homology of Banach algebras, we refer the reader
to \cite{H}.
Recall that a Banach algebra $A$ is amenable if for every Banach $A$-bimodule $E$,
the first order bounded Hochschild cohomology group of $A$ with coefficients
in the dual Banach bimodule $E^*$ vanishes, $\mH^1(A,E^*)=0$, or equivalently
any bounded derivation $D:A\to E^*$ is inner. A famous Theorem proved by
B.E. Johnson \cite{J}, says that for any locally compact group $G$,
amenability of $G$ (\cite{P}) is equivalent to the amenability
of the convolution group algebra $L^1(G)$. For a modern account
on amenability see \cite{Ru}.
\begin{proposition}\label{p1}
Let $A$ and $B$ be Morita equivalent self-induced Banach algebras. Suppose that $A$ is amenable.
Then for every two sided $B$-induced module $E\in B$-$\mathbf{mod}$-$B$, and $n\geq1$, $\mH_n(B,E)=0$,
and the complete quotient seminorm of $\mH_0(B,E)$ is a norm.
\end{proposition}
\begin{proof}
Corollary IV.10 of \cite{G2}.
\end{proof}
\begin{theorem}\label{t5}
Let $G$ be an amenable discrete group, $I$ be a nonempty set and $A=\e(B(I,G))$.
Then for any two sided $A$-induced Banach $A$-bimodule E and every $n>0$, the topological
Hochschild homology groups $\mH_n(A,E)$ are trivial and $\mH_0(A,E)$ is a Banach space.
\end{theorem}
\begin{proof}
It was proved in the preceding section that $\e(B(I,G))\approx\e(G)\oplus\mathbb{C}$.
By Johnson's Theorem, $\e(G)$ is an amenable Banach algebra and thus so is $\e(G)\oplus\mathbb{C}$.
Now, apply Proposition \ref{p1}.
\end{proof}
Note that for any self-induced Banach algebra $A$, the class of two sided $A$-induced modules
is very wide, since for any Banach $A$-bimodule $F$, the module $A\ho_AF\ho_AA$ is two sided $A$-induced,
and also, any $A$-induced bimodule is of this form.

The following Theorem directly follows from duality between definitions of Hochschild homology and cohomology,
Theorem \ref{t5}, and general results of homology theory in the
category of Banach spaces and continuous linear maps, see
for instance \cite{H} or Theorem 4.8 of \cite{GLW}.
\begin{theorem}
Let $G$ be an amenable discrete group, $I$ be a nonempty set and $A=\e(B(I,G))$.
Then for any two sided $A$-induced Banach $A$-bimodule $E$ and every $n>0$, the bounded
Hochschild cohomology groups $\mH^n(A,E^*)$ are trivial.
\end{theorem}
A Banach algebra $A$ is called {\it approximately amenable} \cite{GL}, \cite{GLZ}, if for any Banach $A$-bimodule $E$,
every bounded derivation $D:A\to E^*$ is {\it approximately inner}, that is for a net $(f_\lambda)\in E^*$ and every
$a\in A$, $D(a)=\text{ norm-}\lim_\lambda a\d f_\lambda-f_\lambda\d a$. The following Theorem (that corrects! some
preceding results on amenability properties of Brandt semigroup algebras) is proved in \cite{MP}.
\begin{theorem}\label{t3}
Let $S=B(I,G)$ be a  Brandt semigroup.  Then the following are
equivalent.
\begin{enumerate}
\item [(1)]  $\e(S)$ is amenable.
\item [(2)]  $\e(S)$ is approximately amenable.
\item [(3)]  $I$ is finite and $G$ is amenable.
\end{enumerate}
\end{theorem}
\begin{theorem}
The notion of approximate amenability of self-induced Banach algebras is not a Morita invariant.
\end{theorem}
\begin{proof}
Let $G$ be an amenable group, $I$ be a finite nonempty set and $J$ be an infinite set. Then by Theorem  \ref{t3},
$\e(B(I,G))$ is approximately amenable and $\e(B(J,G))$ is not approximately amenable. But by Theorem \ref{t4},
we have $\e(B(I,G))\approx\e(B(J,G))$.
\end{proof}
\bibliographystyle{amsplain}

\end{document}